\documentclass[a4paper,12pt]{article}
\usepackage[cp1251]{inputenc}
\usepackage{amssymb,amsmath,amsthm}
\usepackage[russian]{babel}
\usepackage{amssymb,amsmath,amsthm}
\usepackage{indentfirst}
\usepackage{hyperref,hypbmsec}
\usepackage[pdftex]{graphicx}
\numberwithin{equation}{section}
\newtheorem{Th}{\hskip\parindent Теорема}[section]
\newtheorem{Le}{\hskip\parindent Лемма}[section]

\newtheorem{Zam}{\hskip\parindent Замечание}[section]
\newtheorem{Hyp}{\hskip\parindent Гипотеза}[section]

\newcommand{\A}{\mathcal{A}}

\newcommand{\E}{\mathfrak{C}}
\newcommand{\R}{\mathfrak{R}}

\newcommand{\D}{\mathfrak{D}}
\newcommand{\N}{\mathbb{N}}
\newcommand{\M}{\mathfrak{M}}
\newcommand{\NN}{\mathfrak{N}}
\newcommand{\1}{\mathbf{1}}

\newcommand{\q}{\mathbf{q}}

\newcommand{\p}{\mathbf{p}}

\newcommand{\Z}{\mathbb{Z}}
\newcommand{\rr}{\mathbb{R}}
\newcounter{propet}

\renewcommand{\le}{\leqslant}\renewcommand{\ge}{\geqslant}

\begin{document}
\author{И.\, Д.\, Кан 
(I.\,D.\,Kan)\footnote{Работа выполнена при поддержке РФФИ (грант  15-01-05700 А)} \\
Московский Авиационный институт \\
(Национальный исследовательский университет)
}

\title{
\begin{flushleft}
УДК 511.321 + 511.31
\end{flushleft} Усиление \, \, теоремы Бургейна --- Конторовича \,-\, V\,\, 
}
\date{}
\maketitle

\date{}
\maketitle
\begin{abstract} Zaremba's conjecture (1971) states that every positive integer number $d$ can be represented as a denominator (continuant) of a finite continued fraction  $\frac{b}{d}=[d_1,d_2,\ldots,d_{k}],$ with all partial quotients $d_1,d_2,\ldots,d_{k}$ being bounded by an absolute constant $A.$ Recently (in 2011) several new theorems concerning this conjecture were proved by Bourgain and Kontorovich. The easiest of them states that the set of numbers satisfying Zaremba's conjecture with $A=50$ has positive proportion in $\N.$ In 2014 Kan and Frolenkov proved this result with $A=5.$ 
 Let $\E_{\A}$  be the set of infinite continued fractions whose partial quotients belong to $\A$
$$\E_{\A}=\left\{[d_1,\ldots,d_j,\ldots]: d_j\in\A,\,j=1,\ldots\right\}$$
and let $\delta$ be the Hausdorff dimension of $\E_{\A}.$ Naw this result  proved  with $A=4$ and 
$\delta>0.7807\ldots$. 

В настоящей работе доказывается, что знаменатели тех конечных цепных дробей, все неполные частные которых принадлежат произвольному конечному  алфавиту $\A$ с параметром $\delta>0.7807\ldots$ (т. е. такому, что множество  бесконечных цепных дробей с неполными частными из этого алфавита
имеет хаусдорфову размерность $\delta$, удовлетворяющую неравенству $\delta>0.7807\ldots$), содержат положительную долю натуральных чисел. Ранее аналогичная теорема была известна лишь для алфавитов с несколько большими значениями  $\delta$. Именно, впервые результат такого рода для произвольного конечного алфавита с $\delta>0.9839\ldots$
 получили в 2011 году Бургейн и Конторович. Далее, в 2013 году автор статьи совместно с Д. А. Фроленковым доказали теорему для  произвольного конечного алфавита с $\delta>0.8333\ldots$. Результат автора 2015 года, предшествующий настоящему, относился к  произвольному конечному алфавиту с $\delta>0.7862\ldots$.

  \noindent
\textbf
{Библиография:} 14 названий.\\
\noindent
\textbf
{Ключевые слова и выражения:} цепная дробь,  тригонометрическая сумма, гипотеза Зарембы, континуант.\par

\end{abstract}

\setcounter{Zam}0
\section{Введение}

\textbf{1.1. История вопроса}

Для $\ d_1,d_2,\ldots,d_{k}\in \N$ 
через $\left     [d_1,d_2,\ldots,d_{k}\right]$ обозначена конечная цепная дробь

\begin{equation}
\label{cont.fraction}
[d_1,d_2,\ldots,d_k]=\cfrac{1}{d_1+\cfrac{1}{d_2+{\atop\ddots\,\displaystyle{+\cfrac{1}{d_k}}}}},
\end{equation}
 а через $\R_{\A}$ --- множество рациональных чисел $\frac{b}{d}$ , представимых  конечными цепными дробями  с {неполными частными}
 $\ d_1,d_2,\ldots,d_{k}$ из некоторого конечного алфавита $\A \subseteq \N$:

  $$
\R_{\A}=\left\{\frac{b}{d}=[d_1,d_2,\ldots,d_{k}]\Bigl|{} \  d_j\in\A{}\ \mbox{для}{}\  j=1,2,\ldots,k\right\}.
$$
Через $\D_{\A}$    обозначено множество несократимых знаменателей $d$  чисел $\frac{b}{d}\in\R_{\A}$: 
$$\D_{\A}=\left\{d\in \N\Bigl|{} \  \exists b\in \N:{}\  \gcd(b,d)=1, {}\ b<d, {}\  \frac{b}{d}\in\R_{\A}\right\}
.
$$
\begin{Hyp}\label{h1.1} (Заремба,  \cite{Zaremba}). Существует константа $A$, такая что  для алфавита 
$\A= \{1,2,\ldots,A\} 
$
имеет место равенство 
$\D_{\A}=\N$.\end{Hyp}

Обзор результатов, связанных с гипотезой 1.1, можно найти в  работах \cite{BK},\cite{NG}. 

Для алфавита $\A$ число $d$ называется \textbf{допустимым} \cite{BK}, если  для любого $q>1$  множество $\D_{\A}$ содержит хотя бы одно число,  сравнимое с $d$ по модулю $q$. Множество допустимых чисел обозначено через $\mathfrak{A}_{\A}$.
Через  $\delta=\delta_{\A}$ обозначим хаусдорфову размерность множества бесконечных цепных дробей с неполными частными из алфавита ${\A}$. Для каждого   элемента $d\in {\D}_{\A}$ его
\textbf{кратностью}
называется количество натуральных чисел чисел $b<d,$ взаимно простых с $d$ и
 таких, что $\frac{b}{d}\in\R_{\A}$. 
Бургейн и Конторович в 2011 году доказали следующее.

\begin{Th}\label{t1.1} \cite[см. теоремы 1.2 и 1.8]{BK}. Пусть 
алфавит ${\A}$ удовлетворяет условию 
$ \delta>\frac{307}{312}=0.9839\ldots.                                                                                                                                   
$
Тогда  множество $ \D_{\A}$ содержит почти все допустимые числа. Точнее, найдется константа $c=c({\A})>0$, такая что для всех достаточно больших чисел  $N$ во множестве  $ \D_{\A}\cap\left [0.5N,N\right]$ содержится по крайней мере 
 \begin{equation}
 \label{lld}                                                                                                                                                                                                             
{\left|\mathfrak{A}_{\A}\cap\left [0.5N,N\right] \right|}
 \left(1-
 \exp{\left\{-c\sqrt{\log N}\right\}}
 \right) 
 \end{equation} 
 элементов,  кратность которых --- не меньше величины
  \begin{equation}
 \label
 {cdd}                                                                                                                                                                                                             
\gg{}\ N^ {2\delta-{1.001}}.                                                                                                                                   
\end{equation}
В частности, справедливо неравенство  
\begin{equation}
 \label{lc}                                                                                                                                                                                                                                                                                                                                                                                                                    
 \left|\D_{\A}\cap\left [1,N\right] \right|\gg N.                                                                                                                                   
\end{equation}                                                                                                              
\end{Th}

   Теорема \ref{t1.1}  применима к алфавиту $\A=\{1,2,\ldots,50\}$ \cite{BK}. 
 В дальнейшем  по поводу различных обобщений теоремы \ref{t1.1} несколькими авторами был написан целый ряд работ  (\cite{FK1} --- \cite{Magg}). 
        Так, в  \cite{FK4}  автор настоящей статьи совместно с Д. А. Фроленковым доказали неравенство (\ref{lc}) для алфавита $\A=\{1,2,3,4,5\}$.          
     
     Далее, Хуанг \cite{Huang} доказал для того же алфавита $\A=\{1,2,3,4,5\}$ формулы  (\ref{lld}) и (\ref{cdd}). 
     Из доказательства теоремы Хуанга следует общий принцип: для всякого $c$, такого что $0.5<c<1$, для вывода   при условии $\delta>c$
     формул (\ref{lld}) и (\ref{cdd})
     методом Бургейна --- Конторовича
 достаточно тем же методом   доказать формулу (\ref{lc}) при том же условии.      
    С помощью этого принципа в  \cite{K6} было доказано, что в теореме \ref{t1.1}  оценку   $ \delta>0.9839\ldots                                                                                                                                   
$
 можно заменить   условием $ \delta>
 \frac{1}{3}\left(\sqrt{19}-2\right)=0.7862\ldots,  
$ которому удовлетворяет алфавит  $\A=\{1,2,3,4\}$.     
    
 Недавно в статье  \cite{Magg}  Магее, Ох и Винтер доказали, что для некоторого положительного числа $\varepsilon$  для алфавита $\A=\{1,2,3,4,5\}$ имеет место неравенство $\left|\D_{\A}\cap\left [1,N\right] \right|\ge 
 N\left(1-
 N^{-\varepsilon}
 \right)$.

\textbf{1.2 Основной результат работы}

Сформулируем основной результат  настоящей статьи. 
\begin{Th}\label{2.1}  
Пусть  алфавит ${\A}$ удовлетворяет условию   \begin{equation}
 \label{hy}                                                                                                                                                                                                             \delta>
 0.25\left(\sqrt{17}-1\right)=0.7807\ldots.                                                                                                                                   
 \end{equation}
 Тогда множество  $\D_{\A}$ содержит положительную долю натуральных чисел и почти все допустимые числа. Точнее, справедливы формулы   
 (\ref{lld})                                                                                                                                                                                                                                                                                                                                                                                                                                                                                                                                                                                                                  --- 
 (\ref{lc}).                                                                                                                                                                                                                                                                                                                                                                                                                    
\end{Th}

          Автор благодарит профессора Н. Г. Мощевитина за постановку задачи и неоднократное обсуждение темы статьи. Также автор благодарен Д. А. Фроленкову за многократное обсуждение и многие  полезные советы.

\section{Основные свойства ансамбля }\label{6}
Через  $G_{\A}\subseteq SL\left(2,\mathbb{Z}\right)$ обозначена мультипликативная полугруппа с единицей  $E=
\begin{pmatrix}
1 & 0 \\
0 & 1
\end{pmatrix}$, состоящая из произведений $2\times2$-матриц:
\begin{equation}
\label{1yrorm}
\begin{pmatrix}
a & b \\
c & d
\end{pmatrix}=
\begin{pmatrix}
1 & d_2 \\
d_1 & d_1d_2+1
\end{pmatrix}
\begin{pmatrix}
1 & d_4 \\
d_3 & d_3d_4+1
\end{pmatrix}\ldots
\begin{pmatrix}
1 & d_k \\
d_{k-1} & d_{k-1}d_{k}+1
\end{pmatrix},
\end{equation}
где $k$ --- четно и $d_1,d_2,\ldots,d_{k}\in\A$. За норму матрицы (\ref{1yrorm}) обычно (\cite{FK1} --- \cite{K5}) принимается   величина
$||g||=d=\langle d_1,d_2,\ldots,d_{k}\rangle
$
 ---  знаменатель конечной цепной дроби (\ref{cont.fraction}), не сократимый с ее числителем (или \textbf{континуант}).

Скажем, что при $n\in\N$ для некоторого множества $\Xi \subseteq G_{\A}$ имеет место разложение  
\begin{equation}
\label{1morm}
 \Xi=\Omega_1\Omega_2\Omega_3\ldots\Omega_n
 \end{equation}
на \textbf{независимые} множители   $ \Omega_1,\Omega_2,\Omega_3,\ldots,\Omega_n\subseteq G_{\A}$, если для каждой матрицы $\gamma\in\Xi$ найдется, причем единственный, набор матриц $g_1,g_2,g_3,\ldots,g_n,$ таких что   
$$\gamma=g_1g_2g_3\ldots g_n,{}\ {}\ g_i\in \Omega_i{}\ {}\ 
 (i=1,2,3,\ldots ,n).
$$
Конечно, при этом выполняется равенство  
$ \left|\Xi\right|=|\Omega_1||\Omega_2||\Omega_3|\ldots|\Omega_n|.
 $

Всюду далее будем использовать следующие обозначения:  $N$   --- достаточно большое натуральное число,  $\varepsilon_0~\in~\left  (0,{}\ 0.0004\right)$ --- фиксированное произвольно малое,  $A=\max\A$, 
$Q_1= \left[\exp{ \left  (A^4\varepsilon_0^{-5}\right)}\right]+1,
 $ 
где для всякого $w\in \rr$ положено $[w]=\min
\{
z\in\Z
\bigl|{}\ 
z\le w
\}$.
 Для каждого $j\in \Z$
 положим 
\begin{equation}
\label{norm9yy} 
Q_j=\1_{\left\{j\not=0\right\}}\left (Q_1\right)^j,\end{equation}
  где, здесь и далее, ${\1_{\{P\}}}=1$ для истинных утверждений $P$ 
   и ${\1_{\{P\}}}=0$ --- для ложных. 

 Рассмотрим  две  произвольные  матрицы
$  g_2\in\Omega_2,{}\  {}\  g_4\in\Omega_4$, три действительных положительных параметра $M_1,M^{(2)},M^{(4)}$
 и  следующие два неравенства: 
\begin{equation}
\label{6unorm}
 \frac{ M^{(2)} }{150A^2{\left  (M_1\right)}^{2\varepsilon_0}} \le||g_2||\le 73A^2   M^{(2)}{ \left  (M_1M^{(2)}\right)}^{2\varepsilon_0},{}\ {}\ {}\ {}\ {}\ 
\frac{{ \left  (M^{(4)}\right)}^{1-\varepsilon_0} }{150A^2} \le||g_4||\le 73A^2M^{(4)}.
 \end{equation}
По достаточно большому числу $N$ и по малому параметру $\varepsilon_0~\in~\left  (0,{}\ 0.0004\right)$  в   \cite{FK3} было построено специальное множество матриц --- \textbf{ансамбль}  (см. терминологию в \cite{BK}) $$\Omega^{(N)}=\Omega^{(N,\varepsilon_0)}\subseteq \left\{\gamma\in G_{\A}\,\Bigl|\,||\gamma||\le1,02 N\right\},$$
 для которого имеет место разложение на независимые множители  (\ref{1morm}) с $n=4$ со свойствами, перечисленными в следующей лемме.            
 
  \begin{Le}\label{l6.1}
 \cite[теорема 3.1]{K5} Существует  непустое множество матриц --- ансамбль $\Omega^{(N)}\subseteq G_{\A}$, такое что для  всякого числа $M_1 \in  \left[Q_1,Q_{-8}N\right]$ найдeтся разложениe 
 $\Omega^{(N)}=\Omega_1\Omega
 $  на независимые множители  $\Omega_1$ и $\Omega$, для которых выполнен ряд свойств: 
   
    во-первых, имеет место оценка 
\begin{equation}
\label{2norm} 
 \left   |\Omega_1\right|\gg{ \left  (M_1\right)}^{2\delta-\varepsilon_0},
  \end{equation}
 
 во-вторых, для любых двух матриц $g_1\in\Omega_1$ и 
$g\in\Omega$  выполняются неравенства 
 \begin{equation}
\label{normyy}
\frac{M_1}{70A^2} \le||g_1||\le 1.01(M_1)^{1+2\varepsilon_0},{}\ {}\ {}\ {}\ 
\frac{N}{160A^2(M_1)^{1+2\varepsilon_0}} \le||g||\le 73A^2\frac{N}{M_1},
 \end{equation}

    в-третьих, для любых чисел $
  M^{(2)},M^{(4)}\in  [Q_1,N]\bigcup\{1\}  $,  таких что
\begin{equation}
\label{yorm}
M_1M^{(2)}M^{(4)}\le N,
  \end{equation}
  найдется разложение $\Omega=\Omega_2\Omega_3\Omega_4$ 
   на независимые множители   $\Omega_2, \Omega_3, \Omega_4$,  для которых  выполнены как неравенства
   \begin{equation}
\label{2n1orm}
 \left   |\Omega_2\right|\gg
 { \left  (M^{(2)}\right)}^{2\delta} 
 { \left  (M_1M^{(2)}\right)}^{-2\varepsilon_0},{}\ 
 \left  |\Omega_4\right |\gg{ \left  (M^{(4)}\right)}^{2\delta-2\varepsilon_0},
  \end{equation}  
так и оценки (\ref{6unorm})  --- для любых двух  матриц $  g_2\in\Omega_2,{}\  {}\  g_4\in\Omega_4$.
 В частности, если $M^{(2)}=1$ или $M^{(4)}~=~1,$ то $\Omega_2=\{E\}$ или $\Omega_4=\{E\}$, соответственно.
\end{Le}

Пусть число  $M_1 \in \left[Q_1,\frac{1}{Q_8}N\right]$ уже  выбрано, так что имеет место разложение $\Omega^{(N)}=\Omega_1\Omega $ 
со свойствами (\ref{2norm}) и (\ref{normyy}). В этом случае множество $\Omega$  будем называть \textbf{полуансамблем}.
Если в $G_{\A}$ имеют место разложения $\Omega=\Omega_2\Omega_{3,4}=
\Omega_2\Omega_{3}\Omega_{4}
$ на независимые множители,
то для любых двух элементов $g'$ и $g$  полуансамбля $ \Omega$ введем обозначения 
$$g'=
g'_2g'_{3,4}
=g'_2g'_3g'_4
,{}\ {}\ {}\ {}\ {}\ {}\ 
g=g_2g_{3,4}
=g_2g_3g_4,
$$
  где нижние индексы  указывают на принадледность соответствующим множествам из разложений.  Далее,  обозначим: если $B$ --- некоторое множество $2\times2$-матриц $\gamma$, то $\widetilde{B}$ --- множество  вектор-столбцов $\widetilde{\gamma}=\gamma\left(0\atop{1}  \right)
$.        
 Для координат трех пар  произвольных векторов
  $$\widetilde{g}',\widetilde{g}\in \widetilde{\Omega},{}\ {}\ {}\ {}\ 
  \widetilde{g}'_{3,4}, {}\ \widetilde{g}_{3,4}\in \widetilde{\Omega}_{3,4},{}\ {}\ {}\ {}\ 
  \widetilde{g}'_4,\widetilde{g}_4\in \widetilde{\Omega}_4
  $$  
   введем такие обозначения:
\begin{equation}
\label{11norm}
 \widetilde{g}'=\begin{pmatrix}
 x \\
X
\end{pmatrix},{}\ 
\widetilde{g}=\begin{pmatrix}
y \\
Y
\end{pmatrix},{}\ {}\ 
\widetilde{g}'_{3,4}=\begin{pmatrix}
 x_{7}  \\
X_{7}
\end{pmatrix},{}\ 
\widetilde{g}_{3,4}=\begin{pmatrix}
y_{7} \\
Y_{7}
\end{pmatrix},{}\ {}\ 
\widetilde{g}'_4=\begin{pmatrix}
x_4 \\
X_4
\end{pmatrix},{}\ {}\ 
\widetilde{g}_4=\begin{pmatrix}
y_4 \\
Y_4
\end{pmatrix}.
\end{equation}
  
  Пусть $M_2,M_4\in \mathbb{R}_+$. Рассмотрим неравенства 
  \begin{equation}
\label{6unorm1}
 { M_2
 { \left  (M_1M_2\right)}^{-6\varepsilon_0} }
 \left  ({Q_5}\right)^{-1} 
 \le ||g_2||\le 
 {M_2},
 {}\ {}\ {}\ {}\ {}\ {}\ 
{{ \left  (M_4\right)}^{1-\varepsilon_0} }\left  ({Q_5}\right)^{-1} 
  \le||g_4||\le {M_4}.
 \end{equation}  
 Свойства полуансамбля $\Omega$ из леммы \ref{l6.1} несколько уточняет  следующая

 \begin{Th}\label{t6.1}
  Для любого числа $M_1
 \in \left[Q_1,Q_{-8}N\right]$ найдeтся   непустое множество матриц --- полуансамбль $\Omega\subseteq G_{\A}$, такое что   для произвольного числа $
  M_2\ge 1   $,  удовлетворяющего неравенству
$M_1M_2\le Q_{-5}N,
 $ найдется разложение  вида $\Omega=\Omega_2 \Omega_{3,4}$   
   на независимые множители   $\Omega_2$ и $\Omega_{3,4}$,  для которых выполнен ряд свойств: 
   
    (i) имеет место  неравенство
   \begin{equation}
\label{2n1orm1}
 \left   |\Omega_2\right|\gg
 { \left  (M_2\right)}^{2\delta} 
 { \left  (M_1M_2\right)}^{-10\varepsilon_0},
  \end{equation}  
  
(ii) для любой  матрицы $g_2\in\Omega_2$ имеет место первая из  оценок (\ref{6unorm1}),
  
    (iii) для любого числа $
  M_4  \ge 1  $,  удовлетворяющего неравенству
\begin{equation}
\label{1yorm1}
M_1M_2M_4\le Q_6 N,
  \end{equation}
    найдется разложение  $\Omega_{3,4}=\Omega_3 \Omega_4$   
   на независимые множители   $\Omega_3$ и $\Omega_4$,  такие что  выполнены  как  неравенство
  \begin{equation}
\label{2n1orm2}
\left  |\Omega_4\right |\gg{ \left  (M_4\right)}^{2\delta-2\varepsilon_0},
    \end{equation}
так  и вторая из оценок в (\ref{6unorm1})  --- для любой  матрицы $g_4\in\Omega_4$: в частности,  
 \begin{equation}
\label{normyy1}
{ \left  (M_4\right)}^{1-\varepsilon_0}  Q_{-5}
\ll 
X_4 \le {M_4}
,{}\ {}\ {}\ {}\ 
{ \left  (M_4\right)}^{1-\varepsilon_0}  Q_{-5}
\ll 
Y_4\le {M_4}
.
 \end{equation}
 \end{Th} 
 Доказательство. 
 Пусть выбраны значения величин $M_2$ и    $M_4$, удовлетворяющих условиям теоремы.  Тогда  положим:
$$\mathcal{M}^{(2)}={M_2}Q_{-3}\left(M_1M_2\right)^{-4\varepsilon_0},{}\ {}\ {}\ {}\ 
\mathcal{M}^{(4)}={M_4}Q_{-3}.
$$
Участвующие в лемме \ref{l6.1} величины $M^{(2)}$ и $M^{(4)}$ определим правилами:

$
M^{(2)}=
                {\mathcal{M}}^{(2)},$ если  ${\mathcal{M}}^{(2)}\ge Q_1,$                
   и           $M^{(2)}=  1$ --- в противном случае;  аналогично,

             $M^{(4)}=
               \mathcal{M}^{(4)},$ если  $\mathcal{M}^{(4)}\ge Q_1,$
         и    $M^{(4)}=  1$ --- в противном случае.
          
    \noindent      Тогда числа  $M^{(2)}$ и $M^{(4)}$ принадлежат множеству $ [Q_1,N]\bigcup\{1\}$. Далее,  если ${M}^{(2)}=~1$ или ${M}^{(4)}=1$, то, полагая $\Omega_2=\{E\}$ или $\Omega_4=\{E\}$, соответственно, получаем, что все  условия   леммы \ref{l6.1} выполнены. Если же  ${M}^{(2)}\not=1$ и ${M}^{(4)}\not=1$, то, следовательно, выполнены равенства $
M^{(2)}=
                {\mathcal{M}}^{(2)}$ и $
M^{(4)}=
                {\mathcal{M}}^{(4)}$. Для таких значений 
   условие  (\ref{yorm})   выполнено ввиду неравенства (\ref{1yorm1}). Поэтому доказаны  оценки (\ref{6unorm}) и (\ref{2n1orm}).
Подставляя в них   значения ${M}^{(2)}$ и ${M}^{(4)}$, получаем неравенства (\ref{6unorm1}), (\ref{2n1orm1}) и (\ref{2n1orm2}). Теорема доказана.

Отметим, что далее довольно часто в качестве значений ${M}_2$  и ${M}_4$ выбираются числа  
\begin{equation}
\label{17narm}
  {M}_2=1,{}\ {}\ {}\ {}\     M_4=\max\left\{1, {}\ {}\  0.4Q_{\alpha-1}\right\} ,
 \end{equation}
где $\alpha\in \N$. В этом случае для проверки оценки (\ref{1yorm1}) достаточно установить неравенство
\begin{equation}
\label{17narb}
  {M}_1 Q_{\alpha} \le 
 {Q_6N}
 .
 \end{equation}

\section{Основа доказательства формул 
 (\ref{lld})                                                                                                                                                                                                                                                                                                                                                                                                                                                                                                                                                                                                                  --- 
 (\ref{lc}).                                                                                                                                                                                                                                                                                                                                                                                                                 
}

\label{7}
Напомним обозначения из \cite{K5}. Применяя теорему Дирихле \cite[лемма 2.1, стр. 17]{Von}, для каждого $\Theta\in[0,1)$ найдем числа $\lambda\in \left (-\frac{1}{4},\frac{1}{4}\right ],$ $q\in\N$ и целые числа $a$  и $l$, такие что 
\begin{equation}
\label{17norm}
  \Theta=\left\{\cfrac{a}{q}+\cfrac{l}{2N}+\cfrac{\lambda}{N}\right\}
  ,{}\ {}\ {}\  \gcd(a,q)=1, {}\ {}\ 
  \frac{q-1}{q}
  \le a<q\le\frac{\sqrt{N}}{Q_1},{}\ {}\ 
   |l|\le\frac{3}{q}Q_1\sqrt{N},
 \end{equation}
где $\{w\}=w-[w]{}\ $ ---  дробная доля числа $w \in \rr.$
 Для аналогичного разложения другого произвольного числа $\Theta'\in[0,1)$
    введем  обозначение:
$\Theta'=\left\{\cfrac{a'}{q'}+\cfrac{l'}{2N}+\cfrac{\lambda'}{N}\right\}.
$ Фиксируем параметр $\lambda$ (в дальнейшем всегда будет $\lambda'=\lambda$), константу
 $T_1=7Q_7
$  и целое число $\kappa\in\left[0,{}\ T_1-1\right]$. Для натуральных индексов  $\alpha$ и $\beta$  положим   
\begin{equation} 
\label{67mn4rm}
P_{\alpha,\beta}=P^{(\lambda)}_{\alpha,\beta}(\kappa)=\left \{\Theta{}\ 
 \mbox{из}{}\  (\ref{17norm}){}\ \Bigl|{} \   l\equiv\kappa\pmod {T_1},{}\ {}\ {}\ 
Q_ {\alpha-1}< q\le Q_{\alpha} , {}\ {}\ 
Q_ {\beta-1}\le         \left|l \right|         \le Q_{\beta}  \right\}.
\end{equation}

Предположим, что число $M_1$ удовлетворяет неравенству
\begin{equation} 
\label{69n1fo2r3m}
 75A^2 Q_{\alpha}Q_{\beta}{}\ {}\ \le {}\ {}\ M_1{}\ {}\ \le {}\ {}\ 
  \min\left\{\left(  Q_{\alpha}Q_{\beta}   \right)^5   ,{}\ 
  \left(Q_{\alpha}Q_{\beta} \right)^{-0.5}{N}
\right\}.
\end{equation}
 Рассмотрим разложение $\Omega^{(N)}=\Omega_{1}\Omega,$ соответствующее этому значению  $M_1$.    Всюду далее   $Z$ --- любое непустое подмножество конечного множества $P_{\alpha,\beta}$. Для любых двух чисел $\Theta',\Theta\in Z$ положим:
\begin{equation} 
\label{67mnorm}
\p=\gcd\left(q',q\right),{}\ {}\ {}\ {}\ q'_0=\frac{q'}{\p},{}\ {}\ {}\ q_0=\frac{q}{\p},{}\ {}\ {}\ 
\q=\frac{q'q}{\p}=\p q'_0q_0,{}\ {}\ {}\  \mathbf{P}= \frac{74A^2Q^2_{\alpha}Q_{\beta}}{ M_1},{}\ {}\ 
 \mathbf{T}=  \frac{\mathbf{P}}{\p}    .
\end{equation}
Используя  обозначения (\ref{11norm}), через $t$ и  $T$ обозначим целые числа, для которых  выполнены соотношения
\begin{equation} 
\label{69n2f1orm}
a'{q_0}x-aq'_0y
\equiv t \pmod{\q}, {}\ {}\ {}\  a'{q_0}X-aq'_0Y\equiv T
\pmod{\q},{}\ {}\ {}\ {}\ {}\  |t|,|T|{}\ \le\frac{\q}{2}.
\end{equation}
Тогда имеют место равенства
${}\ {}\ {}\ 
\left\| \cfrac{a' x}{q'}-\cfrac{ay}{q}\right\|=
\cfrac{|t|}{\q},{}\ {}\ {}\ {}\ 
\left\|  \cfrac{a'X}{q'}-\cfrac{aY}{q}\right\|=
\cfrac{|T|}{\q},$ где  через  
   $\|w\|=\min\left\{\{w\},\{-w\}\right\}$ обозначено расстояние от $w\in \rr$ до ближайшего целого.

Рассмотрим соотношения      
\begin{equation} 
\label{52n4orm}
xT-Xt\equiv 0  \pmod{{q'_0}}, {}\ {}\ {}\ {}\ 
yT-Yt\equiv 0  \pmod{{q_0}},
 \end{equation}
\begin{equation} 
\label{52norm}
\max\{|t|,{|T|}\}\le
74A^2Q_{\beta}
\left({M_1} \right)^{-1}{\q}
\le
\mathbf{T}
\le
\mathbf{P},
 \end{equation}
\begin{equation} 
\label{53norm}
\left| xl'-yl \right|\le\left(9A \right) ^5x+2N{|t|}
{\q} ^{-1}
,{}\ {}\ {}\ {}\ {}\ {}\ {}\ {}\ 
\left| Xl'-Yl \right|\le\left(9A \right) ^5X+2N{|T|}{\q} ^{-1}
,
\end{equation}
\begin{equation} 
\label{52n487m}
\1_{\left\{Q_{\alpha}\ll_{\varepsilon} N^{0.5-\varepsilon}
,{}\ {}\ |t|+|T|\not=0\right\}}
\1_{\left\{xT=Xt\right\}}=0
,{}\ {}\ {}\ {}\ 
\1_{\left\{Q_{\alpha}\ll_{\varepsilon} N^{0.5-\varepsilon}
,{}\ {}\ |t|+|T|\not=0\right\}}
\1_{\left\{yT=Yt\right\}}=0,
 \end{equation}
\begin{equation}  
\label{2487m}
\1_{\left\{Q_{\alpha}\ll N^{\frac{1-\delta}{2\delta-1}}
,{}\ {}\ \delta>0.75,{}\ {}\ |t|+|T|\not=0\right\}}
\1_{\left\{xT=Xt\right\}}=0
.
 \end{equation}
   
 Определим множество $$\NN=
\left\{    \left(\widetilde{g}',\widetilde{g},\Theta',\Theta\right) \in   \left(\widetilde{\Omega}, \widetilde{\Omega},Z,Z\right) 
\Bigl|{} \  \mbox{выполнены  соотношения} {}\ \mbox{(\ref{67mnorm})}{}\  \mbox{---}{}\ (\ref{2487m})\right\}.
$$ 
Для каждого $\Theta \in Z$ положим: 
 $$\NN \left(\Theta\right)=
\left\{    \left(
\widetilde{g}',\widetilde{g},\Theta'
\right) \in   \left(\widetilde{\Omega}, \widetilde{\Omega},Z\right) 
\Bigl|{} \  \left(\widetilde{g}',\widetilde{g},\Theta',\Theta\right)  \in \NN
\right\},
$$ 
$$\NN^{\left(\Theta\right)}_{\p,t,T}=\left\{    \left(\widetilde{g}',\widetilde{g},\Theta'
\right) \in   \NN\left(\Theta\right)
\Bigl|{} \  \mbox{значения параметров}{}\ \p,t,T{}\  \mbox{фиксированы} \right\},$$
\begin{equation} 
\label{55no1rm}
\M\left(\Theta\right)=\left\{    \left(\widetilde{g}',\widetilde{g},\Theta'\right) \in  \NN \left(\Theta\right)  \Bigl|{} \  {}\   T=t=0 \right\},
{} \  {}\  {} \  {}\    {} \  {}\ 
\M_0\left(\Theta\right)=
\sum\limits_{\p \bigl| q} 
\sum_{|t|,|T|{}\    \le  \mathbf{T},
\atop{|t|+|T|\not=0}}
\left|\NN^{\left(\Theta\right)}_{\p,t,T}\right|,
\end{equation}
где сумма по $\p$ распространена на все делители числа $q.$

Напомним, что для всякого $\gamma\in\Omega^{(N)}$ 
положено $\widetilde{\gamma}=\gamma \left(0\atop{1}
\right).$
\begin{Le}\label{l7.3}  Если выполнено неравенство (\ref{69n1fo2r3m}),
то имеет место оценка 
\begin{equation} 
\label{55norm}
\sum_{\Theta\in Z}\left|
\sum_{\gamma\in\Omega^{(N)}}e\left((0,1)
\widetilde{\gamma}\Theta\right)
\right|\ll 
\left(M_1\right)^{1+2\varepsilon_0}   
\sqrt{\left|\Omega_1\right|
\left| \NN   \right| }
.
\end{equation}

\end{Le}
Доказательство. В целом данная лемма была доказана ранее (см. \cite[лемму 5.3]{K5}). Некоторого внимания заслуживают разве что соотношения (\ref{52n4orm}),   (\ref{52n487m}) и (\ref{2487m}), которых  прежде не было. Докажем, что эти соотношения следуют из остальных формул, определяющих  множество $ \NN.$

  Для начала следует 
вычесть  из первого  из  сравнений (\ref{69n2f1orm}), умноженного на $T$, второе из них, умноженное на на $t$.  Тогда получим сравнение 
  \begin{equation} 
\label{69n2f1o6rm}
a'{q_0}(xT-Xt)\equiv
 a{q'_0}(yT-Yt)
  \pmod{\p q'_0
q_0}{}\ \equiv 0 \pmod{ q'_0
}
.
\end{equation}
Но  числа
 $q'_0$ и $a'$  взаимно просты: это следует из несократимости дроби $\frac{a'}{q'}$. Кроме того, числа  $q'_0$ и $q_0$ взаимно просты по построению. Отсюда и из сравнения (\ref{69n2f1o6rm}) получаем    первое из сравнений в (\ref{52n4orm}). Второе --- полностью аналогично.
 
 Далее,  заметим, что 
 в  (\ref{52n487m}) утверждается следующее: если для числа $\alpha$ из (\ref{67mn4rm})  выполнена оценка 
\begin{equation} 
\label{2n4o4r8m}
Q_{\alpha} \ll_{\varepsilon} N^{0.5-\varepsilon}, 
  \end{equation}
 то  при  $ |t|+|T|\not=0$ выполнены неравенства
$tX\not=Tx,{}\ {}\ tY\not=Ty .
$ Докажем это свойство, предположив противное.  Действительно: пусть, скажем, $tX=Tx$. Это означает, что при $t,T\not=0$ выполняется равенство $\frac{x}{X}=\frac{t}{T}.$ Отсюда, ввиду несократимости первой из этих дробей,  из неравенств (\ref{52norm}) и (\ref{2n4o4r8m}) 
  следует цепочка оценок  
\begin{equation} 
\label{52n9o4rm}
 X
 \le 
 T
 \le
  \mathbf{P}
  \ll 
 \frac{Q_{\alpha}^2Q_{\beta}}{M_1}
 =
\frac{Q_{\alpha}(Q_{\alpha}Q_{\beta})}{M_1}
 \ll
\frac{N^{0.5-4\varepsilon_0}
(Q_{\alpha}Q_{\beta})} {M_1}
 \end{equation}
при $\varepsilon=4\varepsilon_0$.
Ввиду неравенств из (\ref{17norm}) и (\ref{67mn4rm}),
 выполнена оценка $Q_{\alpha}Q_{\beta} \ll \sqrt{N}.$ Следовательно, продолжая цепочку неравенств (\ref{52n9o4rm}) и используя  нижнюю оценку, взятую из первого из неравенств в (\ref{normyy1}), получаем:
 $$
 X
\ll
\frac{N^{0.5-4\varepsilon_0}(Q_{\alpha}Q_{\beta})} {M_1}
 \ll
\frac{N^{1-4\varepsilon_0}} {M_1}
\ll 
X
\left(\frac{M_1}{N^2}\right)^{2\varepsilon_0}
\ll 
\frac{X}{N^{2\varepsilon_0}}.
$$ Этим противоречием первое из  равенств в  (\ref{52n487m}) доказано. Второе --- аналогично. 

В частности, если $\delta>0.75,$ то  выполнено неравенство $ N^{\frac{1-\delta}{2\delta-1}}\ll_{\varepsilon} N^{0.5-\varepsilon}.$ Поэтому, если первое из  равенств в (\ref{52n487m}) выполнено, то из него следует равенство  (\ref{2487m}).
Лемма доказана.

\begin{Zam} Легко показать, что имеет место равенство
\begin{equation} 
\label{79n4fo2rm1}
\left|\NN\left(\Theta\right)\right|= 
\left|\M\left(\Theta\right)\right|+\M_0\left(\Theta\right).
\end{equation}
 Действительно, для этого достаточно разбить множество $\NN\left(\Theta\right)$ на составляющие его подмножества и выделить слагаемое с $t=T=0.$  Тогда, ввиду соотношений   (\ref{52norm}) и     (\ref{55no1rm}),
 получаем:  
$$
\left|\NN\left(\Theta\right)\right|= 
\sum_{ \p \bigl| q} 
\left|\NN^{\left(\Theta\right)}_{\p,0,0}\right|
+
\sum_{ \p \bigl| q} 
\sum_{|t|,|T|    \le   \mathbf{T},\atop{|t|+|T|\not=0}}
\left|
\NN^{\left(\Theta\right)}_{\p,t,T}
\right|
=
\left|\M\left(\Theta\right)\right|
+
\M_0\left(\Theta\right).
$$
\end{Zam}   

 Через    $\mathbf{c}=\mathbf{c}(\A)>0$ обозначим  произвольную достаточно малую  константу, зависящую только от   $\A$. 
    С помощью метода Хуанга \cite{Huang} (обобщившего методы  Бургейна --- Конторовича  \cite{BK} и других \cite{FK4}) в \cite{K5} была доказана следующая

\begin{Le}\label{l7.1}(\cite[теорема 4.1]{K5}) Пусть для любых   $\alpha,\beta\in \N$ выполняется оценка  
\begin{equation} 
\label{45norm}
\sum_{\Theta\in Z}\left|
\sum_{\gamma\in\Omega^{(N)}}e\left((0,1)\widetilde{\gamma}\Theta\right)
\right|\ll
   \frac{\left   |\Omega^{(N)}\right|\sqrt{|Z|}}
  {\left(Q_{\alpha}Q_{\beta}\right)^{\mathbf{c}+O(\varepsilon_0)}}.
\end{equation}
 Тогда для алфавита $\A$ имеют место формулы  (\ref{lld})                                                                                                                                                                                                                                                                                                                                                                                                                                                                                                                                                                                                                  --- 
 (\ref{lc}).                                                                                                                                                                                                                                                                                                                                                                                                                 
 \end{Le} 

    Положим
    
 \begin{equation} 
\label{70okr2m11}
\Lambda\left(M_1
 \right)=
\left(M_1
 \right)^{-O_+(\mathbf{c})+O(\varepsilon_0)}.\end{equation}

\begin{Th}\label{T7.1} Пусть для любых  натуральных значений $\alpha$ и $\beta$ найдется число $M_1=M_1(\alpha,\beta),$ удовлетворяющее  неравенствам  (\ref{69n1fo2r3m}) и   
 \begin{equation} 
\label{70okrm11}
\max\limits_{\Theta\in Z}{}\ 
{\left|\M\left(\Theta\right)\right|+\max\limits_{\Theta\in Z}{}\ 
\M_0\left(\Theta\right)}
\ll 
|\Omega|^2\left(M_1
 \right)^{-2+2\delta}\Lambda\left(M_1
 \right).
 \end{equation}
Тогда  для алфавита $\A$ имеют место  формулы
(\ref{lld})                                                                                                                                                                                                                                                                                                                                                                                                                                                                                                                                                                                                                  --- 
 (\ref{lc}).                                                                                                                                                                                                                                                                                                                                                                                                                 
\end{Th}

Доказательство. С помощью соотношений (\ref{79n4fo2rm1}) и   (\ref{70okrm11})   легко получить неравенство 
\begin{equation} 
\label{70nokrm}
{\left| \NN\left (\Theta\right)\right|}{}
\ll 
|\Omega|^2\left(M_1
 \right)^{-2+2\delta}\Lambda(M_1)=
 {\left|\Omega^{(N)}\right|^2}{\left|\Omega_1\right|^{-2}}\left(M_1
 \right)^{-2+2\delta}\Lambda(M_1)
 \ll 
  {\left|\Omega^{(N)}\right|^2}{\left|\Omega_1\right|^{-1}}\left(M_1
 \right)^{-2}\Lambda(M_1)
  \end{equation}
 ввиду оценки (\ref{2norm}). Домножим итог  неравенства  (\ref{70nokrm}) на ${\left|\Omega_1\right|}\left(M_1
 \right)^{2+4\varepsilon_0}.$ Далее, суммируя по $\Theta$ и извлекая корень, получаем:
\begin{equation} 
\label{70nok5r}
\left(M_1\right)^{1+2\varepsilon_0}   
\sqrt{\left|\Omega_1\right|
\left| \NN    \right| }
=
\left(M_1\right)^{1+2\varepsilon_0}   
\sqrt{\left|\Omega_1\right|
\sum_{\Theta\in Z}
\left| \NN (\Theta)   \right| }\ll \left|\Omega^{(N)}\right|
\sqrt{|Z|}\Lambda(M_1)
.
 \end{equation}
 Из неравенств  (\ref{55norm})  и  (\ref{70nok5r})
легко получить  оценку (\ref{45norm}). Поэтому утверждение теоремы следует из леммы \ref{l7.1}. Теорема доказана.
\label{8/1}

\section{Оценка величины $\left|\M\left(\Theta\right)\right|$}
Напомним обозначения  (\ref{11norm}) и   (\ref{17norm}) и положим $r=\gcd(y_{7},q),$ $\mathfrak{u}=\frac{y_{7}}{r}$.   Введем обозначение $\mathfrak{u}^{-1}$ 
 для вычета по модулю $\mathfrak{q}=\frac{q}{r}$, обратного к   $\mathfrak{u}$. 
 Через $g_3{\Omega}_4$ обозначим множество матриц, получающихся умножением матрицы $g_3$ на произвольные элементы множества ${\Omega}_4$.  Пусть для целых $m$ и $n$ при $n\not=0$
\begin{equation}
 \label{11hy}                                                                                                                                                                                                                  \delta_n(m)=
\frac{1}{n}\sum\limits^n_{k=1}\exp\left(\frac{2\pi i}{n} km\right)=
\left\{
              \begin{array}{ll}
1,{}\ {}\ \mbox{если}{}\ {}\ m {}\ {}\ \mbox{делится на}{}\ {}\ n,                \\
               0,{}\ {}\ \mbox{--- в противном случае,}                                                                      
              \end{array}
\right.
       \end{equation}
--- $\delta$-символ Коробова (в честь Н. М. Коробова, пропагандировавшего идею использования формулы (\ref{11hy}), см.  \cite{Korobov}).

\begin{Le}\label{9.91} Пусть имеет место разложение на независимые множители $\Omega=\Omega_2\Omega_{3,4}.$ Тогда для произвольного множества $S\subseteq \widetilde{\Omega}_{3,4}$, для произвольного разложения на независимые множители 
$\Omega_{3,4}=\Omega_3\Omega_{4}$
  имеет место оценка 
  \begin{equation} 
\label{10iho11rm}
\sum\limits_{\left(y_{7}\atop{Y_{7}} \right),\left(x_{7}\atop{X_{7}} \right)
\in S}
\delta_{q}
(Y_{7}x_{7}- y_{7}X_{7})
\ll
  \left|\Omega_3\right|
 \sum\limits_{r \bigl| q} {}\ 
 \sum\limits^{\mathfrak{q}}_{k=1}
 \sum\limits_{g_3
\in \Omega_3}
\left(
 \sum\limits_{\left(y_{7}\atop{Y_{7}} \right)
\in \left(S\cap g_3\widetilde{\Omega}_4\right),\atop{\gcd(y_{7},q)=r
}}
\delta_{\mathfrak{q}}\left(Y_{7}\mathfrak{u}^{-1}-k \right)
\right)^2
.
\end{equation}
\end{Le}
Доказательство.  Заметим, что элементы  числовых пар  $(y_{7},Y_{7})$ и  $(x_{7},X_{7})$ взаимно просты как числители и знаменатели цепных дробей, а  для каждого из ненулевых слагаемых   суммы в левой части неравенства (\ref{10iho11rm}) имеет место сравнение $X_{7}y_{7}\equiv x_{7}Y_{7}\pmod q.$ Поэтому выполнено равенство $\gcd(x_{7},q)=\gcd(y_{7},q)=r.$  Далее, полагая
$\mathfrak{v}=\frac{x_{7}}{r},
$
 получаем:
\begin{equation} 
\label{1ho11rm}
\sum\limits_{\left(y_{7}\atop{Y_{7}} \right),\left(x_{7}\atop{X_{7}} \right)
\in S}
\delta_{q}
(Y_{7}x_{7}- y_{7}X_{7})
=
\sum\limits_{r \bigl| q}{}\ 
\sum\limits_{\left(y_{7}\atop{Y_{7}} \right),\left(x_{7}\atop{X_{7}} \right)
\in \widetilde{\Omega}_{3,4,}\atop{\gcd(y_{7},q)=\gcd(x_{7},q)=r
}}
\delta_{\mathfrak{q}}
\left(Y_{7}\mathfrak{v}- \mathfrak{u}X_{7}\right)
\1_{\left\{ \left(x_{7}\atop{X_{7}} \right),\left(y_{7}\atop{Y_{7}} \right)
\in S\right\}}
.
\end{equation}
  Введем обозначение  
                   $\mathfrak{v}^{-1}$
 для вычета по модулю $\mathfrak{q}$, обратного к    $\mathfrak{v}$. 
Тогда, продолжая равенство (\ref{1ho11rm}),
 получаем: 
  \begin{equation} 
\label{19o11r2m}
 \sum\limits_{\left(y_{7}\atop{Y_{7}} \right),\left(x_{7}\atop{X_{7}} \right)
\in S}
\delta_{q}
(Y_{7}x_{7}- y_{7}X_{7})
 =
\sum\limits_{r \bigl| q}
\sum\limits_{\left(y_{7}\atop{Y_{7}} \right),\left(x_{7}\atop{X_{7}} \right)
\in \widetilde{\Omega}_{3,4,}\atop{\gcd(y_{7},q)=\gcd(x_{7},q)=r
}}
\delta_{\mathfrak{q}}
\left(Y_{7}\mathfrak{u}^{-1}- X_{7}\mathfrak{v}^{-1}\right)
\1_{\left\{ \left(x_{7}\atop{X_{7}} \right),\left(y_{7}\atop{Y_{7}} \right)
\in S\right\}}
.
                       \end{equation}    
      Обозначим вычет каждого из чисел $Y_{7}\mathfrak{u}^{-1}$ и $X_{7}\mathfrak{v}^{-1}$ по модулю $\mathfrak{q}$ через $k$. Тогда, вводя сумму по всем таким $k$ и  переходя в правой части равенства  (\ref{19o11r2m}) к квадрату суммы, получаем:                  
  \begin{equation} 
\label{1ho11r2m}
\sum\limits_{\left(y_{7}\atop{Y_{7}} \right),\left(x_{7}\atop{X_{7}} \right)
\in S}
\delta_{q}
(Y_{7}x_{7}- y_{7}X_{7})
=
\sum\limits_{r \bigl| q }{}\  
\sum\limits^\mathfrak{q}_{k=1}{}\  
\left(
 \sum\limits_{\left(y_{7}\atop{Y_{7}} \right)
\in\widetilde{\Omega}_{3,4}
}
\delta_{\mathfrak{q}}\left(Y_{7}\mathfrak{u}^{-1}-k \right) \1_{\left\{ \left(y_{7}\atop{Y_{7}} \right)
\in S,\atop{  \gcd(y_{7},q)=r}\right\}}
 \right)^2.
 \end{equation}                                                                
                   Поскольку $ \left(y_{7}\atop{Y_{7}} \right)=g_3\widetilde{g}_4,$
 то, применив в правой части равенства (\ref{1ho11r2m}) неравенство  Коши ---  Буняковского по переменной $g_3\in \Omega_3$, получим неравенство 
 (\ref{10iho11rm}).
 Лемма доказана.
 
 \begin{Le}\label{9.917} В условиях леммы \ref{9.91} для любого  $S\subseteq \widetilde{\Omega}_{3,4}$ имеет место оценка \begin{equation} 
\label{10ih1o11rm}
\sum\limits_{\left(y_{7}\atop{Y_{7}} \right),\left(x_{7}\atop{X_{7}} \right)
\in S}
\delta_{q}
(Y_{7}x_{7}- y_{7}X_{7})
\ll
{\left   |\Omega_{3,4}\right||S|}
{\left   (Q_{\alpha}\right)^{-\delta+O(\varepsilon_0)}}
.
\end{equation}
\end{Le}

Доказательство.  Для того, чтобы оценка (\ref{10ih1o11rm}) была нетривиальной, необходимо выполнение условия $\alpha>1.$ Раскроем в правой части неравенства (\ref{10iho11rm}) квадрат суммы, далее  производя суммирование по $k$: 
\begin{equation}
\label{197n5orm}
\sum\limits_{\left(y_{7}\atop{Y_{7}} \right),\left(x_{7}\atop{X_{7}} \right)
\in S}
\delta_{q}
(Y_{7}x_{7}- y_{7}X_{7})
\le {\left|\Omega_3\right|} 
 \sum\limits_{g_3
\in \Omega_3} 
\sum\limits_{r \bigl| q} 
 \sum\limits_{\left(y_{7}\atop{Y_{7}} \right),\left(x_{7}\atop{X_{7}} \right)
\in 
\left(S\cap g_3\widetilde{\Omega}_4\right),
\atop{\gcd(y_{7},q)=\gcd(x_{7},q)=r
}}
\delta_{\mathfrak{q}}
(Y_{7}\mathfrak{v}- \mathfrak{u}X_{7})
. 
 \end{equation}
Производя суммирование по  $r$ в  правой части неравенства (\ref{197n5orm}), получаем:
\begin{equation}
\label{1b8n5rm7}
\sum\limits_{\left(y_{7}\atop{Y_{7}} \right),\left(x_{7}\atop{X_{7}} \right)
\in S}
\delta_{q}
(Y_{7}x_{7}- y_{7}X_{7})
\le {\left|\Omega_3\right|} 
\sum\limits_{g_3
\in \Omega_3}
\sum\limits_{\left(y_{7}\atop{Y_{7}} \right),\left(x_{7}\atop{X_{7}} \right)
\in 
\left(S\cap g_3\widetilde{\Omega}_4\right)}
\delta_{q}
\left(Y_{7}x_{7}- y_{7}X_{7}\right) 
. 
\end{equation}
Заметим, что  вектора $\left(x_{7}\atop{X_{7}}\right)$ и $\left(y_{7} \atop{Y_{7}}\right)$  (в обозначениях (\ref{11norm})) получаются из векторов  $\left(x_{4}\atop{X_{4}}\right)$ и $\left(y_{4} \atop{Y_{4}}\right)$ умножением на матрицу $g_3.$ Поскольку $\det g_3=1,$ то такое умножение не изменяет определитель, составленный  из координат двух векторов. Поэтому для   ненулевых слагаемых в  последней из сумм в (\ref{1b8n5rm7}) выполняется  сравнение $X_{4}y_{4}\equiv x_{4}Y_{4}\pmod q$ (вытекающее из сравнения $X_{7}y_{7}\equiv x_{7}Y_{7}\pmod q$).  

С другой стороны, полагая
 $M_4=0.4\sqrt{Q_{\alpha-1}}
$
и подчиняя разложение $\Omega_{3,4}=\Omega_3\Omega_{4}$ этому выбору, 
 согласно теореме \ref{t6.1},  получаем неравенство
 $$\left|X_{4}y_{4}- x_{4}Y_{4}\right|<
 \left(\max\left\{X_{4},{}\  Y_{4}\right\}\right)^2\le \left( M_4\right)^2< Q_{\alpha-1}\le
 q.$$ 
 Следовательно, выполнено равенство $x_4{Y}_4= X_4{y}_4,$ или $\frac{x_4}{X_4}=
\frac{{y}_4}{{Y}_4}$. Ввиду несократимости этих дробей, имеют место  равенства
$x_4=y_4,$  $X_4=Y_4,$ откуда $x_7=y_7,$  $X_7=Y_7.$ Следовательно, неравенство (\ref{1b8n5rm7}) приводит к оценке
$$\sum\limits_{\left(y_{7}\atop{Y_{7}} \right),\left(x_{7}\atop{X_{7}} \right)
\in S}
\delta_{q}
(Y_{7}x_{7}- y_{7}X_{7})
\le 
{\left|\Omega_3\right|} 
\sum\limits_{g_3
\in \Omega_3}
\sum\limits_{\left(y_{7}\atop{Y_{7}} \right),\left(x_{7}\atop{X_{7}} \right)
\in 
\left(S\cap g_3\widetilde{\Omega}_4\right)}
\1_{\left\{x_7=y_7,{}\ X_7=Y_7\right\}}= 
{\left|\Omega_3\right|} 
|S|. 
$$
Отсюда получаем:
$$
\label{1n5rm437}
\sum\limits_{\left(y_{7}\atop{Y_{7}} \right),\left(x_{7}\atop{X_{7}} \right)
\in S}
\delta_{q}
(Y_{7}x_{7}- y_{7}X_{7})
\le 
{\left|\Omega_3\right|}{|S|} =
\frac{\left|\Omega_{3,4}\right||S|}{\left|\Omega_4\right|}{}
\ll 
\frac{\left|\Omega_{3,4}\right||S|}{\left(M_4\right)^{2\delta+O(\varepsilon_0)}}\ll 
\frac{\left|\Omega_{3,4}\right||S|}{\left(Q_{\alpha}\right)^{\delta+O(\varepsilon_0)}}, 
$$
ввиду неравенства (\ref{2n1orm2}). Лемма доказана. 

Напомним, что число $M_4$ было выбрано выше в доказательстве и положим 
$M_2=\max\left\{1, {}\ {}\ 0.4\sqrt{Q_{\beta-1}}\right\}. $ 
Рассмотрим разложение
$\Omega^{(N)}=\Omega_1\Omega_2\Omega_{3,4} $, соответствующего значениям параметров $M_1$ и $M_2$. 

\begin{Le}\label{9.1} Для всякого числа $M_1$  из интервала 
    (\ref{69n1fo2r3m})  
выполнены оценки
\begin{equation} 
\label{109ho11rm}
\max\limits_{\Theta\in Z}
\frac{ 
\left |\M\left(\Theta\right)\right|}{\left   |\Omega_2\right|}
\ll
\sum\limits_{\left(y_{7}\atop{Y_{7}} \right),\left(x_{7}\atop{X_{7}} \right)
\in \widetilde{\Omega}_{3,4}}
\delta_{q}
(Y_{7}x_{7}- y_{7}X_{7})
\ll
\frac{
\left   |\Omega_{3,4}\right|^2}
{\left   (Q_{\alpha}\right)^{\delta+O(\varepsilon_0)}}=
 \frac{|\Omega|^2}{\left  |\Omega_2\right|^2 \left   (Q_{\alpha}\right)^{\delta+O(\varepsilon_0)}}
.
\end{equation}
\end{Le}
Доказательство. 
  Ввиду сравнений в (\ref{69n2f1orm}) с $t=T=0,$ для  элементов множества $\M\left(\Theta\right)$ выполнены сравнения
$ a'q_0x-aq'_0y\equiv 0 \equiv
a'q_0X- aq'_0Y\pmod{\q}
$,
  откуда следует, что 
\begin{equation} 
\label{102n9orm}
   a'q_0x\equiv 0\equiv a'q_0X\pmod{q'_0}
  ,
\end{equation}
\begin{equation} 
\label{102norm}
\left( a'q_0X\right)y \equiv  aq'_0Yy= 
\left( aq'_0y\right)Y\equiv \left( a'q_0x\right)Y\pmod
{q'q_0}.
\end{equation}
Но числа $a'$ и $q'$ взаимно просты ввиду несократимости дроби $\frac{a'}{q'}$, а числа  $q_0$ и $q'$ --- по построению. Следовательно, сокращая начало и конец цепочек сравнений (\ref{102n9orm}) и (\ref{102norm}) на $a'q_0$, получаем сравнения \begin{equation} 
\label{102nor6m}
x\equiv 0 \equiv X\pmod {q'_0} ,{}\ {}\  
Xy\equiv xY\pmod {q'}.
\end{equation}
 Ввиду взаимной простоты чисел $x$ и  $X$,  из первой пары сравнений в  (\ref{102nor6m}) следует, что $q'_0=1$. Аналогично доказывается равенство  $q_0=1$. Следовательно, $q'=q=\q$.
 
  Далее, согласно \cite[лемме 6.1]{K5}, 
 для всех элементов   множества $\M$ имеет место равенство
$g'_2=g_2.
$ Следовательно,  вектора $\left(x \atop{X}\right)$ и $\left(y \atop{Y}\right)$ (в обозначениях (\ref{11norm})) получаются из векторов  $\left(x_{7}\atop{X_{7}}\right)$ и $\left(y_{7} \atop{Y_{7}}\right)$ умножением на матрицу $g_2.$ Поскольку $\det g_2=1,$ то такое умножение не изменяет определитель, составленный  из координат двух векторов. Поэтому, ввиду последнего из сравнений в  (\ref{102nor6m}), выполняется также сравнение $X_{7}y_{7}\equiv x_{7}Y_{7}\pmod q$.  
  
  Подытожим сказанное, оценивая количество элементов во множестве  $\M\left(\Theta\right)$. Для этого  фиксируем матрицу ${g}'_2={g}_2$  одним из $\left   |\Omega_2\right|$ способов. Количество возможностей выбора  элементов $\left(y_{7}\atop{Y_{7}} \right),\left(x_{7}\atop{X_{7}} \right),$ 
ввиду сказанного выше, выражается суммой  из  (\ref{109ho11rm}). Для
  числа $\Theta \in Z,$ данного в виде  (\ref{17norm}),   числа  $a$, $q$ и $l$  определены.   Согласно \cite[доказательству леммы 3.13]{FK3}, число $a'$ определяется по числу $a$  однозначно, исходя из сравнений $ a'x\equiv ay, {}\ 
a'X \equiv aY\pmod{\q}
$: действительно, это следует из взаимной простоты чисел $x$ и $X$. Число  $q'$ равно числу $q$ ввиду сказанного выше.  
 А величина $l'$ определяется по $l$   не более, чем константой способов, согласно следующим из (\ref{53norm})  неравенствам  
$\left| xl'-yl \right|\ll x,{}\ {}\  
\left| Xl'-Yl \right|\ll X.
$ 
Тем самым доказано первое из неравенств в (\ref{109ho11rm}). Второе из них получается применением леммы \ref{9.917} при $S=\widetilde{\Omega}_{3,4}.$
Лемма доказана.

Напомним обозначение (\ref{70okr2m11}) и рассмотрим оценку
\begin{equation}
 \label{70nokrl1}
  \left({Q}_{\alpha}{Q}_{\beta}\right)^{\delta}
   \gg
{\left(M_1
 \right)^{2-2\delta
 }}\Lambda(M_1).
 \end{equation}
 
 \begin{Th}\label{l9.3} Для всякого числа $M_1$,  удовлетворящего неравенствам  
    (\ref{69n1fo2r3m}) и  (\ref{70nokrl1}), 
 имеют место оценки
\begin{equation}
\label{109n1o11rm}
\max\limits_{\Theta\in Z}{}\ 
{\left|\M\left(\Theta\right)\right|}
\ll
{|\Omega|^2}{\left  (Q_{\alpha}Q_{\beta}\right)^{-\delta}}\left (M_1\right)^{O(\varepsilon_0)}
 \ll  
  \left |\Omega\right|^2
\left(M_1
 \right)^{-2+2\delta}\Lambda(M_1).
\end{equation}
\end{Th}

Доказательство. 
Ввиду леммы \ref{9.1} и оценки (\ref{2n1orm1}), имеет место неравенство
$${\left|\M\left(\Theta\right)\right|}|\Omega|^{-2}
 \ll
 {\left  |\Omega_2\right|^{-1} \left   (Q_{\alpha}\right)^{-\delta+O(\varepsilon_0)}}
 \ll 
 {\left  (Q_{\alpha}Q_{\beta}\right)^{-\delta}}
 \left (M_1\right)^{O(\varepsilon_0)}.
$$ Отсюда и из  (\ref{70nokrl1}) следует  вторая   оценка в  (\ref{109n1o11rm}).  Теорема доказана.

\section{Определение и свойства соответственных  чисел} 
 Пусть для  числа
$M_1$ выполнены неравенства (\ref{69n1fo2r3m}), (\ref{70nokrl1}) и
\begin{equation} 
\label{70n6krm11}
\max\limits_{\Theta\in Z}{}\ 
{\M_0\left(\Theta\right)} {}\ {}\ \ll  {}\ {}\ 
{|\Omega|^2}{\left  (Q_{\alpha}Q_{\beta}\right)^{-\delta}} \Lambda(M_1)
 {}\ {}\ \ll  {}\ {}\ 
{|\Omega|^2}\left(M_1
 \right)^{-2+2\delta}\Lambda(M_1).
  \end{equation}
 Тогда число
$M_1$  назовем \textbf{соответственным} 
для   пары  натуральныx чисел $\alpha$ и $\beta$.

\begin{Le}\label{1lL10.1} Если 
   для любых натуральныx чисел $\alpha$ и $\beta$ найдется соответственное  значение $M_1$, 
то для алфавита $\A$ имеют место формулы (\ref{lld})                                                                                                                                                                                                                                                                                                                                                                                                                                                                                                                                                                                                                  --- 
 (\ref{lc}).                                                                                                                                                                                                                                                                                                                                                                                                                 
\end{Le}
Доказательство. Пусть соответственное  значение $M_1$ найдется. Тогда, применяя теорему \ref{l9.3}, получаем, что оценки  (\ref{109n1o11rm}) выполнены. Остается  подставить результат оценок  (\ref{109n1o11rm}) и  (\ref{70n6krm11}) в теорему \ref{T7.1}.
Лемма доказана.

 Пример  соответственных чисел дает следующая 
 
\begin{Le}\label{t9.1} Пусть для числа $M_1\in [Q_1, Q_{-8}N]$ 
 имеет место неравенство (\ref{70nokrl1})  и выполнен хотя бы один из следующих двух наборов соотношений: 
\begin{equation} 
\label{71norm}
M_1=\sqrt{NQ_{\alpha+2}Q_{\beta}},{}\ \end{equation}
\begin{equation} 
\label{71nor5m}
M_1=75 A^2 Q_{2\alpha}Q_{\beta}
,{}\ {}\ {}\ {}\ 
(Q_{\alpha})^{\frac{5}{2}}(Q_{\beta})^{\frac{3}{2}}\le N
.
\end{equation}
Тогда 
число $M_1$ --- соответственное.
\end{Le}

Доказательство. Пусть число $M_1$ определено любым из двух перечисленных способов. Тогда оценка (\ref{69n1fo2r3m}) получается применением неравенства 
$  Q_{\alpha}Q_{\beta}\le Q_3\sqrt{N},
 $ 
 следующего из формул (\ref{17norm}) и (\ref{67mn4rm}).
 
Далее,  для $M_1$ из  (\ref{71nor5m}) равенство $\M_0\left(\Theta\right)=0$
следует из неравенств (\ref{67mnorm}) и  (\ref{52norm}), ввиду которых
$$
\max\{|t|,{|T|}\}
\le
\frac{74A^2Q^2_{\alpha}Q_{\beta}}{ M_1}<1
 .
 $$  
    А для  $M_1$ из  (\ref{71norm}) равенство 
$\M_0\left(\Theta\right)=0$ получается   дословным повторением доказательства \cite[леммы 5.6]{K5}.
  Лемма доказана.

Рассмотрим неравенство   
\begin{equation} 
\label{5.1kr}
N 
\ge 
\left (
\left (Q_{\alpha}\right)^{
 \frac{2\delta-1}{1-\delta}
}
              \left  (Q_{\beta}\right)^{
\frac{2\delta-1}{1-\delta}
}           
\right)^{1-O_+(\mathbf{c})+O(\varepsilon_0)}     
.
\end{equation}   
 Случай невыполнения неравенства (\ref{5.1kr})    рассматривает следующая

 \begin{Le}\label{t9.1h1}
Пусть  для   чисел $\alpha,\beta \in \N$ 
выполнено  неравенство  
\begin{equation} 
\label{59n1okrm65}
N^{1+O_+(\mathbf{c})+O(\varepsilon_0)}   
\le 
 \left (Q_{\alpha}\right)^{
 \frac{2\delta-1}{1-\delta}
}
              \left  (Q_{\beta}\right)^{
\frac{2\delta-1}{1-\delta}
}  
.            
 \end{equation}
Тогда   соответственное значение $M_1$ найдется. 
\end{Le}
Доказательство. Воспользуемся леммой \ref{t9.1} и проверим выполнение неравенства (\ref{70nokrl1}) для числа $M_1$ из (\ref{71norm}). Для этого достаточно проверить, что 
$$\left({N}Q_{\alpha}Q_{\beta}
\right)^{1-\delta}\ll
\left  (Q_{\alpha}Q_{\beta}\right)^{\delta}.
$$
Но именно такое неравенство следует из условия 
(\ref{59n1okrm65}). Лемма доказана.

Таким образом, неравенство (\ref{5.1kr})
всюду далее можно считать выполненным. 

\begin{Le}\label{T8.1} 
Пусть   выполнено неравенство $\delta>\frac{7}{9}$ и пусть для  натуральных  чисел $\alpha$ и $\beta$  имеет место оценка
\begin{equation} 
\label{59n1okrm1}
\frac{\alpha}{\beta}< 
\frac{3\delta-2}{4-5\delta}
. 
\end{equation}  
Тогда найдется соответственное значение $M_1$. 
\end{Le}

Доказательство. Ввиду неравенства  (\ref{5.1kr}), проверка неравенства в  (\ref{71nor5m}) сводится  к установлению оценки

\begin{equation} 
\label{71n5or5m}
(Q_{\alpha})^{\frac{5}{2}}(Q_{\beta})^{\frac{3}{2}}\le  \left  
              (\left (Q_{\alpha}\right)^{\frac{2\delta-1}{1-\delta}}
              \left  (Q_{\beta}\right)^
              {\frac{2\delta-1}{1-\delta}}
              \right)^{1-O_+(\mathbf{c})+O(\varepsilon_0)} .
\end{equation}
Из неравенства $\delta>\frac{7}{9}$  следует, что $\frac{5}{2}<\frac{2\delta-1}{1-\delta}$ и $\frac{3}{2}<\frac{2\delta-1}{1-\delta},$ то есть, что оценка  (\ref{71n5or5m}) выполнена. 
 Следовательно, 
 оценка в  (\ref{71nor5m}) имеет место.  
 
Определим теперь число $M_1$  равенством в (\ref{71nor5m}) и проверим неравенство (\ref{70nokrl1}):
$$ {\left(\left (Q_{\alpha}\right)^{2}Q_{\beta}
 \right)^{2-2\delta}}
 \ll 
\left({Q}_{\alpha}\right)^{\delta}
\left({Q}_{\beta}\right)^{\delta}
 \left({Q}_{\alpha}{Q}_{\beta}\right)^{O(\varepsilon_0)-O_+(\mathbf{c})},
 $$ или, после упрощения,  
$$ \left
(Q_{\alpha}\right)^
{4-5\delta}
\left(Q_{\beta}\right)^
{2-3\delta}
\ll 
 \left({Q}_{\alpha}{Q}_{\beta}\right)^{O(\varepsilon_0)-O_+(\mathbf{c})}
.
 $$ Логарифмируя, получаем  оценку (\ref{59n1okrm1}). Следовательно, оценка (\ref{70nokrl1}) выполнена. Поэтому утверждение настоящей леммы следует из леммы \ref{t9.1}. 
    Лемма доказана.

 Таким образом, при $\delta>0.75$ неравенство 
  \begin{equation} 
\label{59n1krm1}      
\frac{\alpha}{\beta}\ge
\frac{3\delta-2}{4-5\delta}> 1
\end{equation}     
 всюду далее также можно считать выполненным. 
 
 \begin{Th}\label{Tcdd} Пусть  при $\delta>0.75$  соответственное 
   значение $M_1$ найдется
   для любых натуральных чисел $\alpha$ и $\beta$, удовлетворяющих  неравенствам
    (\ref{5.1kr}) и (\ref{59n1krm1}).
 Тогда для алфавита $\A$ имеют место формулы (\ref{lld})                                                                                                                                                                                                                                                                                                                                                                                                                                                                                                                                                                                                                  --- 
 (\ref{lc}).                                                                                                                                                                                                                                                                                                                                                                                                                 
 \end{Th}
 
 Доказательство. Согласно  лемме \ref{1lL10.1}, достаточно найти соответственные
   числа $M_1$ для каждой пары $(\alpha,\beta)$. 
 Однако при невыполнении оценок (\ref{5.1kr}) или (\ref{59n1krm1}) существование соответственного числа доказано в леммах \ref{t9.1h1} или \ref{T8.1}.
 Теорема доказана.
\label{r7r}

\section{Оценка мощности вспомогательных множеств}

Напомним обозначения (\ref{11norm}) и  (\ref{17norm}). Для ближайших трех лемм положим $\Omega_2=$ $\{E\}.$    Через  $S_q\subseteq \widetilde{\Omega}_{3,4}$ обозначим произвольное множество, такое что для любых  двух элементов $\left(y_{7}\atop{Y_{7}} \right),\left(x_{7}\atop{X_{7}} \right)
\in S_q$ 
выполнено сравнение $X_{7}y_{7}\equiv x_{7}Y_{7}\pmod q.$
\begin{Le}\label{LL9.1} Если выполнены неравенства $M_1
\sqrt{ Q_{\alpha}}<N$  и  $\alpha>1,$ то  имеет место оценка 
\begin{equation} 
\label{199iho11rm}
 |S_q|\ll
{\left   |\Omega\right|}
{\left   (Q_{\alpha}\right)^{-\delta+O(\varepsilon_0)}}.
\end{equation}
\end{Le}
Доказательство. 
Используя $\delta$-символ Коробова (\ref{11hy}), согласно  лемме \ref{9.917}, получаем:
 $$ |S_q|^2=\sum\limits_{\left(y_{7}\atop{Y_{7}} \right),\left(x_{7}\atop{X_{7}} \right)
\in S_q}1
=
\sum\limits_{\left(y_{7}\atop{Y_{7}} \right),\left(x_{7}\atop{X_{7}} \right)
\in S_q}
\delta_{q}
\left(Y_{7}x_{7}- y_{7}X_{7}\right) \ll
 \frac{\left   |\Omega_{3,4}\right||S_q|}
{\left   (Q_{\alpha}\right)^{\delta+O(\varepsilon_0)}}=
 \frac{\left   |\Omega\right||S_q|}
{\left   (Q_{\alpha}\right)^{\delta+O(\varepsilon_0)}}.
 $$
 Сокращая последнее неравенство на $|S_q|,$ получаем оценку (\ref{199iho11rm}). Лемма доказана.

Пусть при фиксированном $q'$ числа $w$ и  $W$ ---  два произвольных вычета по модулю $q'.$ Тогда через  $\mathbf{S}=\mathbf{S}_{w,W}$ обозначим число решений системы из двух сравнений
 \begin{equation}
\label{11499r3m1}
x_7a'\equiv w ,{}\  {}\ 
 X_7a'\equiv W\pmod{q'}
\end{equation}
в переменных $a',x_7, X_7,$ таких что $\left( { x}_7\atop{ { X}_7} \right)\in \Omega_{3,4}$ и $0\le a'< q',$ $\gcd(a',q')=1.$
\begin{Le}\label{LL9.421} Пусть выполнены неравенства $M_1
\sqrt{ Q_{\alpha}}<N$  и  $\alpha>1.$ Тогда  для всякого вектора $\left(x_7\atop{X_7} \right)$ найдется не более одного числа $a',$ удовлетворяющего сравнениям (\ref{11499r3m1}), и
   имеет место оценка 
$$\max\limits_{0\le {}\ w,W{}\ <q}{}\ 
 \mathbf{S}_{w,W}
 \ll
{\left   |\Omega\right|}
{\left   (Q_{\alpha}\right)^{-\delta+O(\varepsilon_0)}}.
$$
 \end{Le}
Доказательство.  Из сравнений (\ref{11499r3m1}) и взаимной простоты чисел $a'$ и $q'$  следуют равенства
$$\gcd\left\{x_7,q'\right\}=\gcd\left\{w,q'\right\},{}\ {}\ 
   \gcd\left\{X_7,q'\right\}=\gcd\left\{W,q'\right\}.  
$$
Следовательно, вводя обозначения
$\mathbf{u}=\gcd\left\{x_7,q'\right\}$, $\mathbf{U}=\gcd\left\{X_7,q'\right\},$ из сравнений  (\ref{11499r3m1}) получаем:
 \begin{equation}
\label{113no4m1}
 \frac{x_7a'}{\mathbf{u}}\equiv \frac{w}{\mathbf{u}}\pmod{\frac{q'}{\mathbf{u}}}, {}\ {}\ {}\ {}\ {}\ {}\ 
  \frac{X_7a'}{\mathbf{U}}\equiv  \frac{W}{\mathbf{U}}\pmod{\frac{q'}{\mathbf{U}}}.\end{equation}
Поскольку,  ввиду взаимной простоты чисел $x_7$ и $X_7,$ числа $\mathbf{u}$ и $\mathbf{U}$ также взаимно просты, то 
  сравнениями (\ref{113no4m1}) число $a'$ определяется однозначно (похожая идея была использована в \cite[доказательстве леммы 3.13]{FK3}). 

 Далее, определим  $S_{q'}$ как множество векторов $\left(x_7\atop{X_7} \right),$ для которых найдется какое-либо $ {a}', $  взаимно простое с  $q'$ и  
 удовлетворяющее сравнениям (\ref{11499r3m1}). Пусть $\left( \mathbf{ x}_7\atop{ \mathbf{ X}_7} \right)$ --- также некоторый вектор из множества  $S_{q'}.$
  Это, по определению, означает, что для некоторого $ \mathbf{a}',$ взаимно простого с  $q',$ выполнены сравнения 
 $$\mathbf{x}_7\mathbf{a}'\equiv w ,{}\  {}\ 
 \mathbf{X}_7\mathbf{a}'\equiv W\pmod{q'}.
 $$
 Отсюда и из  (\ref{11499r3m1}) получаем:
 $$x_7{a}'\mathbf{X}_7\mathbf{a}'
 \equiv wW  \equiv  
 X_7{a}'\mathbf{x}_7\mathbf{a}'
 \pmod{q'}. $$
 Сокращая последнее сравнение на взаимно простые с $q'$ числа $ {a}'$ и $ \mathbf{a}',$ получаем  сравнение $x_7\mathbf{X}_7
 \equiv 
 X_7\mathbf{x}_7
 \pmod{q'}. $ Следовательно, применима лемма \ref{LL9.1}, ввиду которой при $\left|S_{q'}\right|=\mathbf{S}$ выполнено требуемое неравенство   (\ref{199iho11rm}).  Лемма доказана.

Пусть $p$ и  $P$ ---  два произвольных вычета по модулю $q.$ Тогда через  $R
\left  (\Theta\right)=R_{p,P}\left  (\Theta\right)\subseteq$ $\subseteq \widetilde{\Omega}_{3,4}$ обозначим произвольное множество двумерных целочисленных векторов, для любого   элемента которого $\left(y_{7}\atop{Y_{7}} \right)
\in R\left  (\Theta\right)$ 
выполнены сравнения $ay_{7}\equiv p,{}\   aY_{7}\equiv P \pmod {q}.$

 \begin{Le}\label{LL9.21} Если выполнены неравенства (\ref{17narb}) и $\alpha>1$, то  имеет место оценка    
$$
\max\limits_{0\le {}\ p,P{}\ <q}{}\ 
  \max\limits_{\Theta\in Z}{}\ 
\left|R_{p,P}\left(\Theta\right)\right|
 \ll 
{\left|\Omega\right|}{\left(Q_{\alpha}\right)^{-2\delta+O(\varepsilon_0)}}
.
 $$\end{Le}
Доказательство.   Рассмотрим разложение полуансамбля $\Omega=\Omega_{3,4}=\Omega_{3}\Omega_{4},$ 
 соответствующее параметру  $M_4=
 0.4{ Q_{\alpha-1}}$, как в (\ref{17narm}).
 Тогда, ввиду этого разложения,  по обеим координатам векторов 
выполнено сравнение
$$   ag_3\left(y_4\atop{Y_4} \right)
    = 
      \left(y_7a\atop{Y_7a} \right)\equiv 
 \left( p\atop{P}\right)
      \pmod{q}.
$$  Пусть  $a^{-1}$ --- вычет по модулю  $q$, обратный к  $a$.
Тогда,  умножая последнее сравнение на число $a^{-1}$, а также  на матрицу, обратную к  матрице  $g_3$, получаем: 
 \begin{equation}
\label{14n1or4m1}
 \left(y_4\atop{Y_4} \right)\equiv  g_3^{-1}
 \left( p\atop{P}\right) 
a^{-1}\pmod{q}.
 \end{equation}

Заметим, что при заданном значении 
$ g_3 $
 правая часть сравнения  (\ref{14n1or4m1}) 
 определена однозначно. Следовательно, определены  остатки компонент вектора
  $\left(y_4\atop{Y_4} \right)$ по модулю  $q$. Однако, ввиду неравенства (\ref{normyy1}), эти компоненты меньше, чем $q$,  поскольку, согласно выбору числа $M_4$ в (\ref{17narm}), выполняется неравенство
  ${M_4} < Q_{\alpha-1}\le q .$ Поэтому  при заданном значении $ g_3$
 вектор $\left(y_4\atop{Y_4} \right),$ а следовательно, и вектор $\left(y_7\atop{Y_7} \right)$  определены однозначно. Однако  матрицу $g_3$ можно выбрать одним из 
   $\left|\Omega_3\right|$ способов. Следовательно, 
   $$  \left|R\left(\Theta\right)\right|
   \le
   {\left|\Omega_3\right|} 
   = 
{\left|\Omega_{3,4}\right|}{\left|\Omega_4\right|^{-1}}
\ll 
{\left|\Omega_{3,4}\right|}{\left(M_4\right)^{O(\varepsilon_0)-2\delta}}\ll 
{\left|\Omega_{3,4}\right|}{\left(Q_{\alpha}\right)^{O(\varepsilon_0)-2\delta}}=
\left|\Omega\right|\left(Q_{\alpha}\right)^{O(\varepsilon_0)-2\delta}
, 
$$
ввиду неравенства (\ref{2n1orm2}). Лемма доказана.

\section{Число решений системы из двух сравнений}
В этом параграфе будет доказана оценка  (\ref{70n6krm11}) при невыполнении условий леммы \ref{t9.1}.
 
 Пусть $m$, $n$, $\p$, $t $ --- данные   целые числа, в том   числе $m$, $n$ и $\p$ не равны  нулю, а $m$ и $n$ взаимно просты. Рассмотрим в целых переменных $x$ и $y$ сравнение
\begin{equation}
\label{109nor1m}
 mx-ny\equiv t \pmod{ mn\p }. 
\end{equation}

\begin{Le}\label{l10.1} Найдутся целые числа $x^{(0)}$,  $y^{(0)}$ и $k$, зависящие только от значений параметров $m$, $n$, $\p$, $t $, такие что  
$0\le k<\p$ и 
для любого решения $(x,y)$ сравнения (\ref{109nor1m}) выполняются сравнения
\begin{equation}
\label{110nor1m}
 x\equiv x^{(0)}+kn \pmod{ n\p },{}\ {}\ y\equiv y^{(0)}+km \pmod{m\p }.
{}
\end{equation}
\end{Le}
Доказательство. 
Пусть $(x,y)$, $(\mathbf{x},\mathbf{y})$ --- какие-нибудь два решения сравнения (\ref{109nor1m}). Тогда подстановка этих значений в исходное сравнение  при последующем вычитании  результатов этих подстановок дает сравнение  
\begin{equation}
\label{111nor1m}
 m(x-\mathbf{x})\equiv n(y-\mathbf{y}) \pmod{ mn\p}.
{}
\end{equation}
Поэтому, ввиду взаимной простоты чисел $m$ и  $n$, выполняются сравнения 
$ x\equiv \mathbf{x}\pmod{n}$, $y\equiv \mathbf{y} \pmod{m}.
$
Следовательно, найдутся числа $k_1$ и $k_2$, такие что  
\begin{equation}
\label{113nor1m}
 x= \mathbf{x}+k_1 n ,{}\ y= \mathbf{y}+k_2  m.
{}
\end{equation}
Подстановка значений (\ref{113nor1m}) в сравнение (\ref{111nor1m}) приводит к сравнению $k_1\equiv k_2\equiv k\pmod{\p}$, и сравнения (\ref{110nor1m}) доказаны.  
Чтобы прийти к неравенству $0\le k<\p$, остается лишь вместо числа $k$ рассмотреть его остаток от деления на  $\p$. Лемма доказана.

   Напомним обозначение  (\ref{70okr2m11}).

 \begin{Le}\label{LL9.2} Если выполнено неравенство (\ref{17narb}) и $\alpha>1$, то   имеет место оценка 
 \begin{equation}
\label{14n2r3m11}
\max\limits_{\Theta\in Z}{}\ 
 \max\limits_{\p \bigl| q} {}\ 
\max\limits_{0\le{}\ t,T{}\    \le    \mathbf{T},\atop{|t|+|T|\not=0}}
\left|\NN^{\left(\Theta\right)}_{\p,t,T}\right|\p^{-2}
  \ll 
{\left|\Omega\right|^2 Q_{\beta}}{\left  (Q_{\alpha}\right)^{L(1-\delta)-2
}}\Lambda(M_1)
 ,
  \end{equation}  
где $L=2$ при $Q_{\alpha}>N^{0.1}$, и $L=3$ --- в противном случае. 
\end{Le}

Доказательство. Положим $\Omega_2=E,$ тогда $\left( { x}\atop{ { X}} \right)=\left( { x}_7\atop{ { X}_7} \right),$ $\left( { y}\atop{ { Y}} \right)=\left( { y}_7\atop{ { Y}_7} \right).$
Сравнения в (\ref{69n2f1orm})   при $n=q'_0$ и  $m=q_0$ превращаются в сравнения вида (\ref{109nor1m}). Так как   $\p n=q'$,  $\p m=q$, то,  согласно лемме \ref{l10.1}, найдутся целые числа $k$ и $K$ в интервале $[0,{} \ \p-1]$ , такие что
\begin{equation}
\label{114or3m1}
x_7a'\equiv x^{(0)}+kq'_0{}\   \pmod{q'}, {}\ {}\   {}\ {}\   {}\ {}\   {}\ {}\ 
y_7a\equiv y^{(0)}+kq_0{}\  {}\ \pmod{q},{}\  {}\ {}\  {}\ 
\end{equation}
 \begin{equation}
\label{114nr3m1}
 X_7a'\equiv X^{(0)}+Kq'_0{}\  \pmod{q'},{}\  {}\ {}\  {}\ 
Y_7a\equiv Y^{(0)}+Kq_0,{}\  \pmod{q}, 
\end{equation}
  где  $X^{(0)}$ и  $Y^{(0)}$ --- константы, аналогичные величинам  $x^{(0)}$  и $y^{(0)}$ из  леммы \ref{l10.1}. Эти константы зависят только от $q,q',t,T.$
 
 Введем обозначения:
 \begin{equation}
\label{1r4n893m1}
  x^{(0)}+kq'_0=w,{}\  {}\ X^{(0)}+Kq'_0=W,{}\  {}\ 
 y^{(0)}+kq_0=p,{}\  {}\ 
 Y^{(0)}+Kq_0=P.
 \end{equation}
Заметим, что при заданных значениях $q',  \Theta,  t,T
 $ параметры     (\ref{1r4n893m1}) зависят только от величин $k$ и $K$. 
 Тогда сравнения (\ref{114or3m1}) и (\ref{114nr3m1}) перепишутся в виде (\ref{11499r3m1}) и 
$$     y_7a\equiv p,{}\  {}\
 Y_7a\equiv P\pmod{q}
.
$$
 Фиксируя две последних величины, получаем значения остатков компонент вектора
  $\left(ay_7\atop{aY_7} \right)$ по модулю  $q$ (при фиксированном значении $a$). Таким образом,  количество векторов $\left(y_7\atop{Y_7} \right)$ можно оценивать по лемме \ref{LL9.21}, в то время как  количество векторов $\left(x_7\atop{X_7}\right)$ и чисел $a'$ --- по лемме \ref{LL9.421}.

 Подытожим сказанное, оценивая количество всех элементов множества $ \NN_{\p,t,T}\left(\Theta\right),$
 для чего имеется два пути. 

   Согласно первому из них,  сначала  выбирается
  число $q'$ одним из не более, чем $Q_{\alpha}$ вариантов,  затем    ---   пара чисел
  $k$ и $K,$  одним из
        $\p^2$ способов. Тем самым, поскольку числа $t,T,\Theta$ заданы, то параметры $w,W,p,P$
    из     (\ref{1r4n893m1}) определены однозначно. Следовательно,  вектора $\left(y_7\atop{Y_7} \right)=\left(y\atop{Y} \right)$ и $\left(x_7\atop{X_7} \right)=\left(x\atop{X} \right)$ могут быть выбраны не более,  чем $\left|R\left(\Theta\right)\right|$ или $
  \mathbf{S}$  способами, соответственно. Наконец,  число $l'$ можно выбрать не более, чем $2Q_{\beta}$ вариантами.  Перемножая все  названные мощности конечных множеств,  ввиду лемм \ref{LL9.421} и \ref{LL9.21}, получаем оценку:
  \begin{equation} 
\label{69n2f12rm}
	 \frac{1}{\p^2}\left|\NN^{\left(\Theta\right)}_{\p,t,T}\right|
 \ll
 Q_{\alpha}
\left|R\left(\Theta\right)\right|
 \mathbf{S}Q_{\beta}
 \ll
       \frac{Q_{\alpha}\left   |\Omega\right|}
{\left   (Q_{\alpha}\right)^{2\delta+O(\varepsilon_0)}} 
\frac{\left|\Omega\right|Q_{\beta}}{\left(Q_{\alpha}\right)^{\delta+O(\varepsilon_0)}}
= \frac{\left|\Omega\right|^2 Q_{\beta}\Lambda(M_1)}{\left  (Q_{\alpha}\right)^{3\delta-1
}}.
\end{equation}

Второй из обсуждаемых путей оценки
отличается от первого только количеством вариантов для выбора величин                 $\left(x\atop{X} \right)$    и      $q'.$   Именно, сначала вектор $\left(x\atop{X} \right)$ выбирается одним из  не более, чем 
        $\left|\Omega\right|$ вариантов. Далее, согласно        
        соотношениям (\ref{2487m}) и (\ref{5.1kr}),
        определитель $xT-Xt$ не равен нулю. С другой стороны, ввиду первого из 
          сравнений  (\ref{52n4orm}), число $xT-Xt$  делится на $q'_0.$ 
Следовательно,       число $q'_0$ может быть выбрано не более, чем $\ll N^{\varepsilon_0}$ способами
      (поскольку количество делителей натурального числа, не превосходящего $N$, не превосходит величины $\ll_{\varepsilon}  N^{\varepsilon},$ см: \cite{Korobov}). При $Q_{\alpha}>N^{0.1}$ выполнена также оценка $ N^{\varepsilon_0}\ll \left(Q_{\alpha}\right)^
      {O(\varepsilon)}\ll \Lambda(M_1)$.
      Умножая выбранное число $q'_0$
  на заданное число    $\p,$ получаем число        $q'.$  Перемножая все количества выбора  названных величин, ввиду леммы  \ref{LL9.21}, получаем:    
\begin{equation} 
\label{69n2f13rm}
	 \frac{1}{\p^2}
	 \left|\NN^{\left(\Theta\right)}_{\p,t,T}\right|
 \ll 
 \left|\Omega\right| \Lambda(M_1) 
 \left|R\left(\Theta\right)\right|
 {Q_{\beta}} =
 \left|\Omega\right|\Lambda(M_1) 
\frac{\left|\Omega\right|}{\left(Q_{\alpha}\right)^{2\delta+O(\varepsilon_0)}}Q_{\beta} \ll
  \frac{\left|\Omega\right|^2 Q_{\beta}\Lambda(M_1)}{\left  (Q_{\alpha}\right)^{2\delta}}.
  \end{equation} 
  Объединяя  неравенства (\ref{69n2f12rm}) и (\ref{69n2f13rm}), получаем оценку  (\ref{14n2r3m11}). 
  Лемма доказана.

\begin{Le}\label{lkdL10.1} Пусть 
   для   натуральных чисел $\alpha$ и $\beta$  найдется число $M_1$ из интервала (\ref{69n1fo2r3m}), такое что
    выполнена оценка (\ref{70nokrl1}) и  хотя бы одно из следующих двух неравенств  
\begin{equation} 
\label{70dokrm11}
\max\limits_{\Theta\in Z}{}\  
 \max\limits_{\p \bigl| q} {}\ \frac{1}{\p^2}
\max\limits_{0\le{}\ t,T{}\    \le    \mathbf{T},\atop{|t|+|T|\not=0}}
 \left|\NN^{\left(\Theta\right)}_{\p,t,T}\right|
 \ll
{
|\Omega|^2}
 {\left(M_1
 \right)^{2\delta}}
 {\left({Q}_{\alpha}\right)^{-4}  
   \left({Q}_{\beta}\right)^{-2}}\Lambda(M_1),
  \end{equation}  
\begin{equation} 
\label{70cokrm11}
\max\limits_{\Theta\in Z}{}\ 
 \max\limits_{\p \bigl| q} {}\ 
\sum\limits
 _{0\le{}\ t,T{}\    \le    \mathbf{T},\atop{|t|+|T|\not=0}}
 \left|\NN^{\left(\Theta\right)}_{\p,t,T}\right|
 \ll
 {|\Omega|^2}
 \left(M_1 \right)^
{2\delta-2 }\Lambda(M_1).
  \end{equation}
 Тогда имеет место оценка (\ref{70n6krm11}), то есть, число $M_1$ --- соответственное. 
\end{Le}
Доказательство.  Рассмотрим величины
$\left|\NN^{\left(\Theta\right)}_{\p,t,T}\right|$
и, пользуясь неравенством  (\ref{70dokrm11}), оценим их сумму по $t$ и $T$ в пределах от  $0$  до   $\mathbf{T}$. Тогда получим оценку (\ref{70cokrm11}). Далее, суммируя неравенства (\ref{70cokrm11}) по всем $\p$, делящим число $q$, получаем оценку (\ref{70n6krm11}).
 Лемма доказана.

\begin{Zam} \label{200fkrm}  Согласно теореме \ref{Tcdd},
для выполнения  неравенства   (\ref{17narb})  и второй из верхних оценок в  (\ref{69n1fo2r3m}) достаточно потребовать, чтобы выполнялись  оценка   (\ref{5.1kr})  и 
неравенства

\begin{equation} 
\label{5g091km1}
M_1\sqrt{Q_{\alpha}Q_{\beta}}\le M_1Q_{\alpha}\le 
\left (
\left (Q_{\alpha}\right)^{\frac{2\delta-1}{1-\delta}}
              \left  (Q_{\beta}\right)^{\frac{2\delta-1}{1-\delta}} \right)^{1-O_+(\mathbf{c})+O(\varepsilon_0)}
             \end{equation}
(первое из которых выполнено  ввиду  оценки  (\ref{59n1krm1})).
\end{Zam}

Рассмотрим равенство 
\begin{equation}
\label{1g4nor1m}
M_1=
\left(
\left(Q_{\alpha}\right)^{\frac{2-\delta}{\delta}}
\left(Q_{\beta}\right)^{\frac{3}{2\delta}}
\right)^{1+O_+(\mathbf{c})+O(\varepsilon_0)}.
\end{equation}
Иллюстрацией к применению рассматриваемого метода служит следующая

\begin{Th} \label{1g4nort} Пусть выполнены неравенства 
$\delta>0.25\left({\sqrt{17}-1}\right)=0.7807\ldots$, $Q_{\beta}<Q_{\alpha}$ и $Q_{\alpha}>N^{0.1}$. Тогда   число  $M_1$ из (\ref{1g4nor1m}) является 
соответственным. 
 \end{Th}

Доказательство. 
Ввиду леммы \ref{lkdL10.1}, достаточно доказать неравенства (\ref{17narb}), (\ref{69n1fo2r3m}),   (\ref{70nokrl1}) и (\ref{70dokrm11}). 
Для этого рассмотрим следующие из условий  (\ref{hy}) и (\ref{59n1krm1})
  неравенства
 \begin{equation}
\label{g4nor1m}
  \left(\frac{2-\delta}{\delta}-
  \frac{\delta}{2-2\delta}\right)\alpha<\left(
  \frac{\delta}{2-2\delta}-\frac{3}{2\delta}\right)\beta,{}\ {}\ {}\ {}\ {}\ {}\ 
   \frac{2}{\delta}<\frac{2\delta-1}{1-\delta},{}\ {}\ {}\ {}\ {}\ {}\ 
\frac{3}{2\delta}<\frac{2\delta-1}{1-\delta}.
 \end{equation}                                                                                      
  Действительно, эти неравенства                                             получаются, соответственно,  из оценок 
    $$\frac{\alpha}{\beta}\ge \frac{3\delta-2}{4-5\delta}>1>\frac{\delta^2+3\delta-3}
    {\delta^2-6\delta+4}  ,{}\ {}\ {}\   
    2\delta^2+\delta-2>0,{}\ {}\ {}\ {}\  
                      4\delta^2+\delta-3>0,$$               
               справедливых, соответственно, при 
                     $$\delta>\frac{7}{9}=0.7777\ldots,{}\ {}\ {}\ {}\                       \delta>\frac{\sqrt{17}-1}{4}=0.7807\ldots,{}\ {}\ {}\ {}\ 
   \delta>\frac{3}{4}=0.7500.$$                   
                                                                                                                                                 
  Подставляя оценки (\ref{g4nor1m}) в показатели степеней, получаем  неравенства   
$$\left(Q_{\alpha}\right)^{\frac{2-\delta}{\delta}}
\left(Q_{\beta}\right)^{\frac{3}{2\delta}}
\le
\left(
\left(Q_{\alpha}Q_{\beta}\right)^{
\frac{\delta}{2-2\delta}}
\right)^{1-O_+(\mathbf{c})+O(\varepsilon_0)},
{}\ {}\ {}\ 
\left(Q_{\alpha}\right)^
{\frac{2}{\delta}}
\left(Q_{\beta}\right)^{\frac{3}{2\delta}}
\le  
\left  
              (\left (Q_{\alpha}Q_{\beta}\right)^\frac{2\delta-1}{1-\delta}
              \right)^{1-O_+(\mathbf{c})+O(\varepsilon_0)} 
          ,
 $$
из которых, далее,  оценки (\ref{70nokrl1}) и  (\ref{5g091km1})
 следуют непосредственно. Поэтому, согласно  замечанию  \ref{200fkrm}, имеют место неравенство (\ref{17narb}) и вторая из верхних оценок в  (\ref{69n1fo2r3m}). Остальные оценки в (\ref{69n1fo2r3m}) получаются применением неравенств
$1<\cfrac{2-\delta}{\delta}<\cfrac{3}{2\delta}<5,
$
выполненных  ввиду оценки $1>\delta>0.25\left({\sqrt{17}-1}\right)$.   
      Следовательно, имеют место условия леммы \ref{LL9.2}, ввиду которой  выполняется   оценка  (\ref{14n2r3m11}). Но для $M_1$ из  (\ref{1g4nor1m})
      эта оценка совпадает с  (\ref{70dokrm11}).    
   Весь список требуемых неравенств доказан.
Теорема доказана.
\label{r7r1}
      
     Таким образом, далее можно считать, что выполнено неравенство $Q_{\alpha}\le N^{0.1}$.        Другими словами, неравенство  (\ref{17narb})  и вторая из верхних оценок в  (\ref{69n1fo2r3m}) теперь выполнены.

\section{Обобщение теоремы \ref{1g4nort} на случай малых $Q_{\alpha}$} 
 Для усиления  результата  теоремы \ref{1g4nort} рассмотрим равенство 
\begin{equation}
\label{1g..r1m}
M_1=
\left(
\left(Q_{\alpha}\right)^{\frac{5-3\delta}{2\delta}}
\left(Q_{\beta}\right)^{\frac{3}{2\delta}}
\right)^{1+O_+(\mathbf{c})+O(\varepsilon_0)}.
\end{equation}

\begin{Th} \label{1g.ort} Пусть выполнены неравенства 
$\delta>0.25\left({\sqrt{17}-1}\right)=0.7807\ldots$ и $Q_{\beta}<Q_{\alpha}\le N^{0.1}$. Тогда   число  $M_1$ из (\ref{1g..r1m}) является 
соответственным. 
 \end{Th}

Доказательство. Ввиду леммы \ref{lkdL10.1}, достаточно доказать неравенства (\ref{17narb}), (\ref{69n1fo2r3m}),   (\ref{70nokrl1}) и (\ref{70dokrm11}). 
Аналогично первому из  неравенств  (\ref{g4nor1m}), рассмотрим неравенство
 \begin{equation}
\label{g.or1m}
  \left(\frac{5-3\delta}{2\delta}-
  \frac{\delta}{2-2\delta}\right)\alpha<\left(
  \frac{\delta}{2-2\delta}-\frac{3}{2\delta}\right)\beta
  .
 \end{equation}                                                                                      
  Действительно, это неравенство                                             получается  из оценки 
    $$\frac{\alpha}{\beta}\ge \frac{3\delta-2}{4-5\delta}>\frac{\delta^2+3\delta-3}
    {\delta^2-6\delta+4},$$               
               справедливой при 
                     $$  \delta>\frac{3+\sqrt{31}}{11}=0.7788\ldots.$$                 
  Подставляя элементы оценки (\ref{g.or1m}) в показатели степеней, получаем  
неравенство   
$$\left(Q_{\alpha}\right)^{\frac{5-3\delta}{2\delta}}
\left(Q_{\beta}\right)^{\frac{3}{2\delta}}
\le
\left(
\left(Q_{\alpha}Q_{\beta}\right)^{
\frac{\delta}{2-2\delta}}
\right)^{1-O_+(\mathbf{c})+O(\varepsilon_0)},
 $$
из которого, далее,  оценка (\ref{70nokrl1}) 
 следует непосредственно. 
 
 Поскольку выполнены неравенства $Q_{\beta}<Q_{\alpha}\le N^{0.1}$, то для доказательства оценки (\ref{17narb}) и второй из верхних оценок в  (\ref{69n1fo2r3m}) достаточно учесть цепочку неравенств 
 $$M_1Q_{\alpha}
 =
 \left(
\left(Q_{\alpha}\right)^{\frac{5-3\delta}{2\delta}}
\left(Q_{\beta}\right)^{\frac{3}{2\delta}}
\right)^{1+O_+(\mathbf{c})+O(\varepsilon_0)}
Q_{\alpha}
 \le 
 \left(Q_{\alpha}\right)^{1.1\left(\frac{5-3\delta}{2\delta}+\frac{3}{2\delta}+1\right)}
 \le 
 \left(N\right)^{0.11\left(\frac{5-3\delta}{2\delta}+\frac{3}{2\delta}+1\right)}<N.
 $$ 
 Остальные оценки в (\ref{69n1fo2r3m}) получаются применением неравенств
$1<\cfrac{5-3\delta}{2\delta}<\cfrac{3}{2\delta}<5,
$
выполненных  ввиду оценки $1>\delta>0.25\left({\sqrt{17}-1}\right)$.   
      Следовательно, имеют место условия леммы \ref{LL9.2}, ввиду которой  выполняется   оценка  (\ref{14n2r3m11}). Но для $M_1$ из  (\ref{1g..r1m})
      эта оценка совпадает с  (\ref{70dokrm11}).    
   Весь список требуемых неравенств доказан. 
  Теорема доказана.

\section{Доказательство  теоремы \ref{2.1}}
Это доказательство следует непосредственно из теорем \ref{Tcdd}, \ref{1g4nort} и \ref{1g.ort}. Действительно, согласно теореме \ref{Tcdd}, достаточно доказать существование соответственных чисел при выполнении условий  (\ref{5.1kr}) и (\ref{59n1krm1}).  Но для больших и малых значений величины $Q_{\alpha}$ существование таких чисел доказано, соответственно, в теоремах  \ref{1g4nort} и \ref{1g.ort}. 
Теорема \ref{2.1} доказана.

\end{document}